\documentclass[a4paper,11pt]{article}
\usepackage{amsmath,amssymb,fullpage,amsthm}
\usepackage{graphicx}
\newcommand{\R}{\mathbf{R}}
\newcommand{\Z}{\mathbf{Z}}
\newcommand{\Q}{\mathbf{Q}}
\newcommand{\N}{\mathbf{N}}

\newtheorem{Theorem}{Theorem}

\newtheorem{Lemma}[Theorem]{Lemma}
\newtheorem{Proposition}[Theorem]{Proposition}

\begin{document}
\title{Non-equilibrium phase diagram for a model with  coalescence, evaporation and deposition}
\author{Colm Connaughton, Centre for Complexity Science, University of Warwick, UK\\
R. Rajesh,  Institute of Mathematical Sciences, CIT Campus, Chennai, India\\
Roger Tribe, Mathematics Institute, University of Warwick, UK\\
Oleg Zaboronski, Mathematics Institute, University of Warwick, UK}
\date{15.11.2012}
\maketitle
\begin{abstract}
We study a $d$-dimensional lattice model of diffusing coalescing massive particles, with two
parameters controlling deposition and evaporation of monomers.
The unique stationary distribution for the system
exhibits a phase transition in all dimensions $d \geq 1$ between a growing phase, in
which the expected mass is infinite at each site, and an exponential phase in which
the expected mass is finite. We establish rigorous upper and lower bounds
on the critical curve describing the phase transition for this system, and some asymptotics
for large or small deposition rates.
\end{abstract}
\vspace{2pc}
\section{Introduction} \label{s1}
In section \ref{s1.1} we state the rules for the interacting particle system studied throughout the paper,
describe the main new results on the phase transition, and explain the intuition behind the results.
In section \ref{s1.2} we describe some previous work on this and related models,
including numerical and theoretical intuition into the two phases.
\subsection{Summary of main results}   \label{s1.1}
\noindent
\textbf{Particle Rules.} Particles live on $\Z^d$. Each particle has a mass with value in $\N = \{1,2,\ldots\}$. Particles
move between nearest neighbour sites as independent rate one simple random walks. When a particle
moves onto an occupied site, the two particles instantly coalesce producing a single particle whose mass
is the sum of the masses of the two interacting particles.  Thus each site contains either zero or one particle at all times.
There are two additional dynamics controlled by parameters $p,q \geq 0$.\\
\vspace{.05in}
\hspace{.1in} \textit{Evaporation.} Rate $p$ evaporation, independently for each particle, reduces their mass by one.
 A particle of mass one therefore disappears at the time of an evaporation.\\
\vspace{.05in}
\hspace{.1in} \textit{Deposition.} At rate $q$, independently at each site $x \in \Z^d$, particles of mass one,
 often called monomers, are deposited. A monomer deposited onto an occupied site
will instantly coalesce, increasing the mass of the existing particle by one.\\
\vspace{.05in}
There is a unique stationary distribution for this process, and the law of any solution
converges to this stationary distribution (see Proposition \ref{basics}).
However the properties of the stationary distribution depend strongly on the values of the parameters
$(p,q)$.
In the next section we summarize some earlier work which describes two phases for this stationary distribution,
called the \textit{exponential phase} and the \textit{growing phase}.
In this paper we characterize these phases by the expected mass at a site $x \in \Z^d$.
Note the stationary distribution is translation invariant.
\begin{Theorem} \label{t1}
Let $M_{\infty}(0)$ have the distribution of the mass of the particle at the origin in the stationary distribution.
Then, there is a non-decreasing function $q_c :[0,\infty) \to (0,\infty) $ so that if $q \in [0,q_c(p))$ then $E[M_{\infty}(0)] < \infty$,
while  if $q > q_c(p)$ then $E[M_{\infty}(0)] = \infty$.
\end{Theorem}
The existence and non-triviality of $q_c(p)$ will be shown, in section \ref{s2}, by analysing the equations
for the first and second moments of the mass at a fixed site.
The moment equations are not closed. However by exploiting global properties of the distribution, namely monotonicity and
pairwise negative correlation, they yield a differential inequality for the first moment that implies the
non-triviality of $q_c(p)$. This was a surprise to the authors,
and it seems to be the presence of negative correlations, natural for systems undergoing coalescence, that allows
the moment equations to be replaced by differential inequalities that are in the correct direction to yield
non-trivial information about the phase curve.
Contrast this with the contact process, which has positive correlations and where moments equations do not seem to help
in establishing the phase transition.
For large or small values of the parameters $p,q$ the system should simplify.  In section \ref{s3} we
give quite detailed theoretical reasoning behind the following asymptotics.

\vspace{.4cm}
\noindent
\textbf{Conjecture.} \textit{Large $p,q$ asympotics.} In all dimensions $q_c(p) < p$ and
\begin{equation} \label{largepq}
 \lim_{p \to \infty} \frac{p-q_c(p)}{p^{1/2}} = \beta_c(d)\in (0,\infty).
\end{equation}
\textit{Small $p,q$ asympotics.} In dimension $d=1$
\begin{equation} \label{smallpq1}
\lim_{p \to 0} q_c(p) p^{-3/2} = \alpha_c \in (0,\infty),
\end{equation}
and in dimensions $d \geq 3$
\begin{equation} \label{smallpq2}
\lim_{p \to 0} q_c(p) p^{-2} = (4 p_d)^{-1}
\end{equation}
where $p_d$ is the escape probability for simple random walk on $\Z^d$, that is the probability
that a walk leaving the origin never returns to the origin.
For each of the three limits, there is a continuum approximation that suggests the answer.
The intuition behind the asymptotic for large $p,q$ is that the many mass changes
occur between each random walk step, so that the model can be approximated,
on suitable time-mass scales, by
a simpler model with continuously evolving masses. We will show that this
simpler model has a one-parameter phase transition with critical value $\beta_c(d)$.
For small $p,q$ in $d=1$ there are many random walk steps between each
deposition and the model can be approximated, on suitable space-time scales, by
one with coalescing massive Brownian particles, and this model too has a  one-parameter
phase transition with critical value $\alpha_c$.
For small $p,q$ in transient dimensions $d \geq 3$, where there is no
continuum approximation for the dynamics. However the when $p$ is small and $q= \alpha p^2$
the occupation density in the stationary distribution is $O(p)$. By rescaling the stationary distribution
in space by $p^{1/d}$ we expect a continuum approximation which has a compound Poisson distribution,
that is particles are positioned at a Poisson rate $\hat{s}_{\alpha}$, and the particles have independent
masses attached with a law determined by a generating function $\hat{\phi}_{\alpha}$.
There is a phase transition for these compound Poisson limits which can be analysed, following the route
used in \cite{MKB1} for the mean field analysis of this model. The critical value for the continuum
approximation can be exactly found as
$\alpha = (4 p_d)^{-1}$, suggesting the asymptotic (\ref{smallpq2}).
The idea that there are various one parameter
models that reflect the limiting situations in the two parameter phase diagram
is similar to the situation for a branching
reacting system studied
in \cite{tribe+M}.
In this paper we do not establish all the weak convergence arguments needed to establish
the conjectured exact asymptotics, but
the moment methods are already sufficient to imply the following weaker forms of the
asymptotics for some cases:
\begin{Theorem} \label{t2}
\begin{eqnarray}
0 <  \liminf_{p \to \infty} \frac{p-q_c(p)}{p^{1/2}}
& \leq & \limsup_{p \to \infty} \frac{p-q_c(p)}{p^{1/2}} < \infty,
\quad \mbox{in all $d$,} \label{crude1} \\
0 <  \liminf_{p \to 0} q_c(p) p^{-3/2}
& \leq &  \limsup_{p \to 0} q_c(p) p^{-3/2} < \infty, \quad \mbox{in $d=1$.} \label{crude2}
\end{eqnarray}
\end{Theorem}
\subsection{Background for model}   \label{s1.2}
An earlier paper on this model \cite{CRO10}
gives a brief bibliography for models of various physical phenomena that involve diffusion, coalescence
(or aggregation) and deposition (or immigration) including the survey \cite{MEA1992}.
The river network model in \cite{scheidegger1967}
and the Takayasu model \cite{privman1997} are close to our model but with zero evaporation. This zero evaporation model
is well understood and has a
(non-reversible) stationary distribution with a constant mass flux from small to large masses.
The mass distribution at a fixed site has polynomial tail, which behaves
like $m^{-4/3}$ for large $m$ in $d=1$ and $m^{-3/2}$ in $d \geq 3$. See \cite{RM2000,CRZtmshort} for recent
developments.
The effect of evaporation was studied in \cite{MKB1}, which demonstrated numerically the phase transition
that we study in this paper and included a mean field analysis.
A more detailed description of the two phases is given in \cite{CRO10},\cite{CRZinout1}
where the mass balance in stationarity is analysed formally, and confirmed by numerical investigations.
In the exponential phase, where $q < q_c(p)$, the mass
distribution at a site has exponential tails. In the growing phase, where $q>q_c(p)$ the mass
distribution at a site has polynomial tail and mimics the tails when there is no evaporation.
Indeed it is conjectured that the entire space-mass distribution,
at large masses, should be well approximated by that of the model with zero evaporation but with
a modified deposition rate $q' = q-p s(\infty)$, where $s(\infty)$ is the occupation probability for
any site in the stationary distribution.
Note that in \cite{CRO10} the phase transition is considered as $\hat{q}_c(p) = \inf \{q: s(\infty) < q/p \} $.
(One can show, see appendix {\ref{s4.1}, that  $s(\infty) < q/p$ for $s > \hat{q}_c(p)$ so that this defines a
true transition).
If $s(\infty) < q/p$ then the first moment equation shows that the first moment is infinite
in the stationary distribution, and hence $q_c(p) \leq \hat{q}_c(p)$.
We do not know whether $q_c(p) = \hat{q}_c(p)$,
and we believe this might be a useful tool in proving some of the many predictions about this model:
exponential moments throughout the exponential phase, behaviour at criticality (see
the scaling predictions in \cite{CRO10}), regularity of the critical curve $p \to q_c(p)$ e.t.c.
\begin{figure}[h]
 \centering
 \includegraphics[scale=1.0]{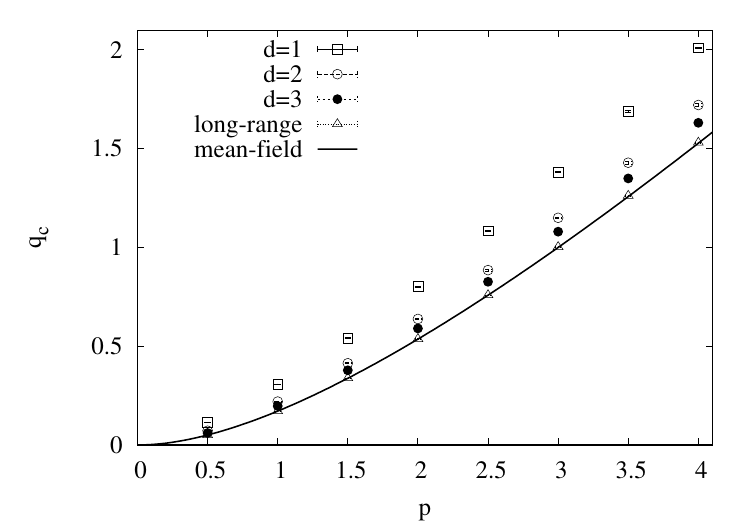}
\caption{\label{fig:phase} Numerical investigation of the phase boundary.
Squares, white circles and black circles are numerical measurements of $q_c(p)$ in $d=1,2,3$ correspondingly.
Triangles are numerical measurements of phase boundary in the approximate '$d=\infty$' model. The solid curve is the rigorous lower bound. }
\end{figure}
Numerical simulations shown in Fig. \ref{fig:phase} illustrate the relation between our bounds and the true phase
boundary in $d=1,2,3$. It is often observed
that mean field approximations become exact in the limit $d \to \infty$ and so one may
conjecture that the critical curve approaches the lower bound (\ref{lowerbound}) as $d$ grows.
An approximate $d=\infty$ model has been simulated using particles that jump
to any other lattice site with equal probability, and the simulated critical curve for this model supports the large
$d$ conjecture.
\section{The phase diagram}  \label{s2}
In section \ref{s2.1} we state monotonicity properties, which follow from standard methods as in
\cite{liggett}, and negative correlation property for the model.
In section \ref{s2.2} we use these and the moment equations to show the existence of the critical curve
$q_c(p)$. The proof of negative correlation, based on the BKR inequality as in \cite{vandenberg+K}, is given in
section \ref{s2.3}.
\subsection{Tools}   \label{s2.1}
We write $M_t(x)$ for the mass of the particle at site $x$ at time $t$.
The construction of the process $ M_t = (M_t(x): x \in \Z^d)$ started from
suitable initial conditions follows from standard arguments on particle models (see for example
chapter IX of \cite{liggett}). One method is to consider $M$ as the solution of the following system of
stochastic differential equations driven by Poisson processes and indexed over the lattice:
almost surely, for all $t \geq 0, \, x \in \Z^d$
\begin{equation}
dM_t(x) =  dP^{(q)}_t(x) - \chi(M_{t-}(x) >0) \, dP^{(p)}_t(x)
 - \sum_{y \sim x} M_{t-}(x) \, dP_t(x,y)
+ \sum_{y \sim x} M_{t-}(y) \, dP_t(y,x)   \label{sde}
\end{equation}
where $y \sim x$ means that $y$ and $x$ are nearest neighbours, and
the independent families of Poisson processes $(P^{(q)}(x): x \in \Z^d)$ of rate $q$ control the deposition,
$(P^{(p)}(x): x \in \Z^d)$ of rate $p$ control evaporation and  $(P(x,y): x,y \in \Z^d)$ of rate $(2d)^{-1}$
control the random walk steps.
We consider only solutions where for all $x$ the path $t \to M_t(x)$ is cadlag
(right continuous with left limits). Furthermore solutions should be adapted to a filtration
$(\mathcal{F}_t)$ where, for any $s<t$, the increment $P_t-P_s$ of any driving Poisson processes
 is independent of $\mathcal{F}_s$.
Such solutions can be found without any further restrictions on the initial conditions. Indeed
we may take as state space the product space $ \mathcal{S} = \N^{\Z^d}$ (with the product topology).
The point is that for independent random walkers without coalescence there would
normally be a growth condition at infinity on the initial state, but this is unnecessary for
instantly coalescing particles (due to the bounded jump intensities).
We collect together some basic results on existence, uniqueness, stationary distributions, moments, monotonicity and correlation, giving comments on the proofs at the end of this section.
\begin{Proposition} \label{basics}
\begin{enumerate}
\item[(i)]
For any initial condition $M_0 \in \mathcal{S}$ that is independent of the driving Poisson processes, there exists a pathwise unique solution to (\ref{sde}) with initial condition $M_0$.
The laws of solutions form a Markov family with a Feller semigroup.
\item[(ii)]
If for some $\theta>0$ and $k \in \N$ the initial moments $E \sum_x |M_0(x)|^k \exp(-\theta |x|)$ are finite,
 then one has finite moments at all times
\begin{equation} \label{moments}
E  \sup_{t \leq T} \sum_x |M_t(x)|^k e^{-\theta |x|} < \infty, \quad \mbox{for all $T$.}
\end{equation}
\item[(iii)]
There exists a unique stationary measure on $\mathcal{S}$. For any solution as above,
the law of $M_t $ converges in distribution, as $t \to \infty$, to the stationary distribution.
\end{enumerate}
\end{Proposition}
A stochastic monotonicity property holds, frequently used for interacting particle systems, as follows.
Put a partial order on $\mathcal{S}$ by writing
$\eta \leq \overline{\eta}$ for if  $\eta(x) \leq \overline{\eta}(x)$
for all $x \in \Z^d$. A function $F:\mathcal{S} \to \R$ is called non-decreasing if $F(\eta) \leq F(\overline{\eta})$
whenever $\eta \leq \overline{\eta}$.
\begin{Lemma} \label{monotonicity}
Let $(M_t:t \geq 0)$ be the solution to (\ref{sde}) started from $M_0 \equiv 0$. Then
for all  measurable and non-decreasing $F: \mathcal{S} \to [0,\infty)$ the expectation
$E\left[ F(M_t) \right]$
is non-decreasing in $t$ and $q$ and non-increasing in $p$.
\end{Lemma}
%
%
\begin{Lemma} \label{negativedependence}
Let $(M_t:t \geq 0)$ be the solution to (\ref{sde}) started from a deterministic initial condition
$M_0 = \eta \in \mathcal{S}$. Then for any $t \geq 0$ the variables
$(M_t(x): x \in \Z^d)$ are negatively associated. In particular,
\begin{equation}
E \left[f(M_t(x))g(M_t(y)) \right] \leq E \left[f(M_t(x)) \right] \, E \left[g(M_t(y))\right]
\label{pairwisenc}
\end{equation}
for measurable and increasing $f,g: \N \to [0,\infty)$.
\end{Lemma}
\noindent
\textbf{Comments. 1.}
Proposition \ref{basics} and the monotonicity Lemma \ref{monotonicity} are rather standard results
for (suitable) particle systems (for example see chapter 2 of Liggett \cite{liggett}). We give some details
in appendix \ref{s4.1}, where we emphasize the use of differential equation comparison methods.
Indeed there is a basic comparison theorem for solutions to (\ref{sde}): if two solutions $(M_t:t \geq 0)$
and $(\overline{M}_t:t \geq 0)$ satisfy $M_0 \leq \overline{M}_0$ a.s., then
$M_t \leq \overline{M}_t$ for all $t \geq 0$ a.s. This can be obtained by standard differential equation methods,
namely a Gronwall estimate
on $E \sum_x e^{-\theta|x|} \chi(M_t(x) > \overline{M}_t(x))$. A similar comparison theorem,
with suitably coupled Poisson drivers, yields the monotonicity in $p,q$ in Lemma \ref{monotonicity}.
The monotonicity in $t$ comes from the fact that $M_0 \equiv 0$ is a minimal initial condition.

\vspace{.05in}
\noindent
\textbf{2.}
The monotonicity can be used to show the existence of a stationary distribution, again in familiar way for
attractive particle systems. Convergence starting from zero initial condition follows from
the monotonicity of  the Laplace transforms
$ E [ \exp(- \sum_{x} M_t(x) \phi(x))]$, for $\phi:\Z^d \to [0,\infty)$ of compact support.
One needs however to check that the limit variables
$(M_{\infty}(x):x \in \Z^d)$ are non-degenerate. Informally, if
$P[M_{\infty}(x) = \infty]>0$ then there would be a steady state
for an extended particle system that allowed particles with infinite mass.
But in this extended system the infinite mass particles do not feel collisions with finite mass particles, and
hence act as an autonomous coalescing system, with no immigration. Such a system does not have
a non-zero steady state.
Pushing this argument further, one can find a maximal entrance law for the extended system
where all sites start with infinite mass. The decay of the infinite mass particles leads to the coupling of the
maximal and minimal entrance laws and this implies the facts about the stationary distribution.
Details are in the appendix.

\vspace{.05in}
\noindent
\textbf{3.}
Negative association is the strongest of various negative correlation type properties
for random vectors, with pairwise negative correlation (\ref{pairwisenc})
being perhaps the weakest.
Pairwise negative correlation, which is all we
need for our arguments,
was established for coalescing random walks in Arratia \cite{arratia} by a Markov duality argument.
We will prove Lemma \ref{negativedependence} by a small modification of the
arguments in van den Berg and Kesten \cite{vandenberg+K}, which deals with purely coalescing
system, and which uses the van den Berg, Kesten, Reimers (BKR)
inequality. Although their argument goes through without much trouble for our model,
we do have to check that the extra deposition and evaporation mechanisms do not destroy the proof.
This is a key tool for this paper, and rather less standard, and so we give the proof in section \ref{s2.3}.
Some restriction on initial conditions is necessary. The examples
in Liggett \cite{liggett2} show that negative correlation properties do not mesh well
with random initial conditions, which is one reason why, unlike positive correlations, differential
techniques have not been that successful in establishing such properties.
\subsection{Phase transition via moments equations}   \label{s2.2}
Throughout this section we consider zero initial conditions. The distribution of $(M_t(x): x \in \Z^d)$ at each fixed $t$ is then
translation invariant on $\Z^d$ (by the uniqueness of solutions to (\ref{sde})).
Monotonicity, as in Lemma \ref{monotonicity}, implies that the moments
$E[M_t^k(0)]$, for $k>0$, are non-decreasing in $t$ and $q$, and non-increasing in $p$.
Set
\[
q_c(p) = \inf \left\{ q: \lim_{t \to \infty} E[M_t(0)] =\infty\right\}
\]
where the expectation uses zero initial conditions and parameter values $(p,q)$.
This defines $q_c(p)$ as a non-decreasing function with values in $[0,\infty]$. Once we have shown that
$q_c(p) \in (0,\infty)$ the remaining statement in Theorem \ref{t1} follows from
monotonicity.
We use the differential equations for the first moment $m_1(t) = E[M_t(0)]$
and the second moment  $m_2(t) = E[M^2_t(0)]$. These can be derived from
(\ref{sde}) by developing $M_t^k(0)$ using calculus and then taking expectation and exploiting
translation invariance.
They imply that the first moment $m_1(t)$ is continuously differentiable and satisfies
\begin{equation} \label{m1}
\frac{dm_1}{dt}(t) = q - p\, s(t)
\end{equation}
where $s(t) = P[M_t(0)>0]$, the probability the origin is occupied at time $t$.
In particular if $q > p$ then $m_1(t) \uparrow \infty$
showing that the steady state will be in the growing phase where $E[M_{\infty}(0)] = \infty$.
This shows $q_c(p) \leq p$, but it is simple to improve this upper bound.
Indeed
\[
\frac{ds}{dt}(t) = q(1-s(t)) - p P[M_t(0)=1] + P[M_t(0)=0, \, M_t(e) > 0] - s(t)
\]
where the last two terms arise from the simple random walking into and out of the origin, and $e$
denotes a nearest neighbour to the origin.
Bounding
\[
P[M_t(0)=0, \, M_t(e) > 0] \leq P[M_t(0)=0] = 1-s(t)
\]
and discarding the term involving $p$ one reaches
\[
\frac{ds}{dt}(t) \leq (1+q) - (2+q) s(t).
\]
Since $s(0) = 0$ this implies that $s(t) \leq (1+q)/(2+q)$ for all $t \geq 0$ and substituting this into
the first moment equation (\ref{m1}) one sees that  $m_1(t) \uparrow \infty$ whenever
$q(2+q) > p (1+q)$. This in turn shows that
\begin{equation} \label{upperbound}
q_c(p) \leq \frac{(p-2) + \sqrt{p^2+4}}{2}.
\end{equation}
Further incremental improvements like this are possible on the upper bound, but without any particular  hope that
one is approaching the true critical curve.
For a lower bound we start with the second moment equation
\begin{equation} \label{m2}
\frac{dm_2}{dt}(t) = q(2m_1(t)+1) - p(2m_1(t)-s(t)) + 2 E[ M_t(0) M_t(e)]
\end{equation}
where the last term arises from the particles leaving or entering the origin.
By monotonicity $dm_2/dt(t) \geq 0$ and by negative correlation $ E[ M_t(0) M_t(e)] \leq m_1^2(t)$.
With these substitutions we find
\[
2m_1^2(t) - 2(p-q) m_1(t) + q +p s(t) \geq 0.
\]
The first moment equation (\ref{m1}) and monotonicity of $m_1(t)$ imply that $s(t) \leq q/p$ for all $t$. This last substitution yields
a quadratic inequality for the first moment:
\begin{equation} \label{qe}
m_1^2(t) -(p-q) m_1(t) + q \geq 0 \quad \mbox{for all $t \geq 0$.}
\end{equation}
Note that the initial condition is $m_1(0) = 0$ and that $t \to m_1(t)$ is continuous and non-decreasing.
Therefore whenever the quadratic in
(\ref{qe}) attains a strictly negative minimum on $(0,\infty)$ the first moment $m_1$ stays bounded by the first positive real root
for all $t \geq 0$.  This happens precisely if $p > q$ and $(p-q)^2 > 4q$, and the corresponding upper bound is
\[
m_1(t) \leq \frac{(p-q) - \sqrt{(p-q)^2-4q}}{2} \quad \mbox{for all $t \geq 0$.}
\]
Resolving the condition $(p-q)^2 > 4q$ with respect to $q<p$ we find that when
$ q < p+2 - 2 \sqrt{p+1}$ the first moment stays bounded and hence
\begin{equation} \label{lowerbound}
q_c(p) \geq p+2 - 2 \sqrt{p+1}.
\end{equation}
Note that this analysis is dimension independent. Moreover the formula for the lower bound appears in the mean
field analysis performed in \cite{MKB1}. This is no surprise since the negative correlation allows us to replace the
second moment equation by a differential inequality that agrees with the mean field approximation.
It is natural to try to apply the same differential inequality methods to investigate other moments, but we have
not profited much from this.  Fro example, exponential moments $E[e^{\theta M_{\infty}(0)}]$ are conjectured to be finite, for suitable $\theta >0$,
whenever $q < q_c(p)$.
Consider the following exponential moment, under the zero initial condition, for a fixed $z \geq 1$.
\[
\psi (z,t) = E [z^{M_t(0)}-1]
\]
This solves formally (since the finiteness of such moments is not yet established) the equation
\begin{eqnarray*}
\frac{d\psi}{dt}(t) &=&  - (pz^{-1}-q)(z-1) \psi(t) + (q-ps(t)z^{-1})(z-1)  \\
& & \hspace{.2in} +  E [z^{M_t(0)+M_t(e)} -   z^{M_t(0)} - z^{M_t(e)} + 1].
\end{eqnarray*}
Monotonicity implies that $\psi(t)$ is increasing and negative correlation shows that
\[
E [z^{M_t(0)+M_t(e)} -   z^{M_t(0)} - z^{M_t(e)} + 1] \leq \psi^2(t).
\]
 Using this we obtain
 the following inequality involving the exponential moment $\psi(t)$ and the
non-occupation probability $s(t)$:
\[
\psi(t)^2 - (pz^{-1}-q)(z-1) \psi(t) + (q-ps(t)z^{-1})(z-1) \geq 0 \quad \mbox{for all $t \geq 0$.}
\]
Now the previous upper bound on $s(t)$ is acting in the wrong direction, but using $s(t) \geq 0$ we reach
the quadratic $\psi(t)^2 - (pz^{-1}-q)(z-1) \psi(t) + q(z-1) \geq 0$. As before, if this quadratic in $\psi$
has a strictly negative minimum at a positive value of $\psi$ then the exponential moment stays bounded.
It is easy to see that this holds when $q>0$ is small, for suitable $z \in (1,q/p)$. This suggests
that the method could be used to show some information on exponential moments, but only
for $q$ far from the critical curve.
\subsection{Proof of negative association}  \label{s2.3}
Various notions of negative correlation are  studied, and they have been exploited for a variety of
interacting particle systems. See Newman \cite{newman}, Pemantle \cite{pemantle}, Liggett \cite{liggett2}
for overviews. Real variables $(X_1,\ldots,X_n) $ are called negatively associated if for
any two disjoint subsets $J_1,J_2 \subseteq \{1,\ldots,n\}$ and any non-decreasing $f_i:\R^{J_i} \to [0,\infty)$ and
one has
\[
E[ f_1(X_j;j \in J_1) f_2(X_j:j \in J_2) ] \leq E[ f_1(X_j;j \in J_1) ] E[ f_2(X_j:j \in J_2) ].
\]
Infinite vectors of variables are called negatively associated if each finite subset is
negatively associated.
We give here the details for the proof of Lemma \ref{negativedependence}, following the strategy from
Lemmas 2.4-2.7 in van den Berg and Kesten \cite{vandenberg+K}. They exploit the BKR inequality,
an inequality for product measures on finite sets, which we state here.
Suppose $V,S$ are finite sets and that $\mu$ is a product measure on the product space $S^V$.
For $K \subset V$ and $\omega = (\omega_v:v \in V) \in S^V$ define
$[\omega]_K = \{ \overline{\omega} \in S^V: \overline{\omega}_v = \omega_v
\; \mbox{for} \; v \in K\}$. The set $[\omega]_K$ is called a cylinder set with base set $K$.
For $A,B \subseteq \Omega$ define
$A \Box B$ to be those $\omega$ for which there exist disjoint $K,K' \subseteq V$ so that $[\omega]_K \subseteq A$
and $[\omega]_{K'} \subseteq B$. The BKR inequality states (see
\cite{reimers}) that $\mu[A \Box B] \leq \mu[A] \, \mu[B]$.
Fix $t>0$. A graphical construction for the process over $[0,t]$ can be given using the same
Poisson processes  that
drive the equation (\ref{sde}). This graphical construction contains, in a direct fashion, information on the
genealogy of individual particles, and we recap the required notation. For each $x \sim y$,
at each jump of $P(x,y)$,  place on the space time lattice $[0,t] \times \Z^d$ an arrow starting
at $(t,x)$ and ending at $(t,y)$. This arrow corresponds to the particle at $(t-,x)$, if it exists,
jumping to site $(t,y)$. Using only a realization of these arrows
there is a natural notion of a path from $(s,x) \to (s',x')$ when $0 \leq s <s' \leq t$. If such a path exists a particle
at $(s,x)$ would end up at $(s',x')$ (ignoring the mass labels for the moment).
See, for example, chapter 3 of Durrett \cite{durrett} for a careful definition.
Using the other Poisson drivers $(P^{(q)}(x): x \in \Z^d)$, $(P^{(p)}(x): x \in \Z^d)$
we add Poisson points where evaporations and depositions occur. Almost surely, one can then trace
through the evolution of the masses $(\eta(x): x \in \Z^d)$ at time zero to yield the
final masses $(M_t(x):x \in \Z^d)$.
We define an embedded discrete time structure which will yield the finite product structure
required for BKR. Fix $N,L \in \N$. Define variables $(N(k,x), Z(k,x): 0 \leq k \leq N-1, x \in \Z^d)$ as follows:
let $I_{k,N}$ be the time interval $ (\frac{kt}{N},\frac{(k+1)t}{N}]$; write $P((r,s])$ for the increment
$P_s-P_r$; let
\[
N(k,x) = \sum_{y \sim x} P(x,y)(I_{k,N})  +  P^{(p)}(x)(I_{k,N}) + P^{(q)}(x)(I_{k,N}),
\]
namely the number of Poisson events that occur at $x$ during the time interval $I_{k,N}$; and let
\[
Z(k,x) = \left\{ \begin{array}{cl}
 0 & \mbox{if $N(k,x)=0$,} \\
p & \mbox{if $N(k,x)=1$ and $P^{(p)}(x)(I_{k,N}) =1$,} \\
q & \mbox{if $N(k,x)=1$ and $P^{(q)}(x)(I_{k,N}) =1$,} \\
\pm e_i & \mbox{if $N(k,x)=1$ and $P(x,x \pm e_i)(I_{k,N})=1$,} \\
 \Delta & \mbox{if $N(k,x) \geq 2$,}
\end{array}
\right.
\]
where $(\pm e_i:i=1,\ldots,d)$ are the unit vectors in $\Z^d$.
Let $S$ be the finite  set of labels $\{0,p,q,\pm e_i,\Delta\}$. Let $V = \{0,\ldots,N-1\} \times (\Z^d \cap [-L,L]^d)$.
Then the vector
\[
\mathcal{Z} := (Z(k,x): 0 \leq k \leq N-1, |x| \leq L)
\]
has the desired product law on $S^V$.
Fix finite disjoint sets $J_1,J_2 \subseteq \Z^d$. Choose $L_0$ so that
$|x| \leq L_0$ for all $x \in J_1 \cup J_2$. Fix  non-decreasing $f_i:\R^{J_i} \to [0,\infty)$,
 and $a_i \geq 0$, for $i=1,2$, and set
\[
\hat{A}_i = \{ f_i(M_t(x_j): j \in J_i) \geq a_i\}, \quad \mbox{for $i=1,2$.}
\]
Negative association follows if we can show $P[\hat{A}_1 \cap \hat{A}_2] \leq P[\hat{A}_1] P[\hat{A}_2]$.
We define a good set $\mathcal{G} = \mathcal{G}^1_{L,N} \cap \mathcal{G}^2_{L}$
as follows.
\[
\mathcal{G}^1_{L,N} = \left\{ N(k,x) + \sum_{y \sim x} N(k,y) \in \{0,1\} \; \mbox{for all $(k,x) \in V$} \right\}
\]
is the set where at most one Poisson event occurs during each interval  $I_{k,N}$ for each $|x| \leq L$, and
furthermore  if one occurs then neighbouring sites have zero Poisson events during this interval.
\[
\mathcal{G}^2_{L} = \left\{ \mbox{there is no path $(x,s) \to (y,t)$ for any $s \in [0,t)$, $|y| \leq L_0$ and
$|x| > L$} \right\}
\]
is the set where no particle moves from outside $[-L,L]^d$ to affect the values
of $(M_t(y): |y| \leq L_0)$.

\vspace{.05in}
\noindent
\textbf{Claim 1.} (i) $P[\mathcal{G}^2_{L}] \to 1$ as $L \to \infty$; (ii) $P[\mathcal{G}^1_{L,N}] \to 1$ as $N \to \infty$.
Write $(X^x_t:t \geq 0)$ for a simple rate one random walk started at $x$. To estimate
$P[\mathcal{G}^2_{L}]$ it is enough,
for each $x$ with $|x| >L$, to consider only paths started at the points $(x,0)$ or $(x,s)$ where
$s <t $ is a jump time of $P(x,y)$.  The expected number of such points that start paths
leading to a $(y,t)$ with $|y| \leq L_0$ is bounded by
\[
\sum_{|x| > L} (1 + t) \sup_{s \leq t}  P[X^x_s \in [-L_0,L_0]^d]
\]
and simple random walk estimates show this approaches zero as $L$ grows.
Part (ii) is a simple Poisson calculation.
\vspace{.05in}
Let $\Omega_0$ be the (good) subset of $S^V$ where the value $\Delta$ is never taken and where
no pair of neighbouring sites, $(k,x), \, (k,y)$ with $x \sim y$, simultaneously take values different from
$0$.  On $\mathcal{G}^1_{L,N}$ we know that $\mathcal{Z} \in \Omega_0$.

\vspace{.05in}
\noindent
\textbf{Claim 2.}
There exist $A_1,A_2 \subseteq \Omega_0 \subseteq S^V$ so that
\[
\hat{A}_i \cap \mathcal{G} =  \{\mathcal{Z} \in A_i\} \cap \mathcal{G} \quad \mbox{for $i=1,2$.}
\]
Indeed, we may take  $A_i = \{\omega \in \Omega_0: \{\mathcal{Z} = \omega\} \cap \mathcal{G} \subseteq
\hat{A}_i \cap \mathcal{G}\}$.
To see this, it is sufficient check that on the set $ \mathcal{G}$ the values of $(M_t(x):|x| \leq L_0)$
can be reconstructed only from $\mathcal{Z}$
and the  the initial masses $(\eta(x): |x| \leq L)$. Of course, this is the purpose of the good set.
For $\omega \in \Omega_0$ we may define (inductively in $k$) mass values $(m(k,x,\omega): |x| \leq L, k=0,\ldots,N-1)$,
starting at $k=0$ with the initial masses $(\eta(x): |x| \leq L)$. This is the natural discrete analogue of our process.
If $\omega_{k,x}$ takes the value $p$ or $q$
the value $m(k-1,x,\omega)$ can be updated to $m(k,x,\omega)$ in the appropriate manner
(since no neighbouring site will affect the mass at $x$ during this interval).
If $\omega_{k,x} = \pm e_i$  the mass $m(k,x,\omega)$ is set to $0$, and the mass $m(k-1,x,\omega)$ is added
to the mass at $y=x \pm e_i$ (if $|y| \leq L$).  Again the restriction to $\Omega_0$  ensures there is never a
conflict due to a pair of events occuring which would require us to decide which
action (deposition, evaporation, movement) to do first.
The mass values $m(k,x,\omega)$ do not take into account any mass entering
from sites $|y| >L$. But on the set $\mathcal{G}^2_L$ we may ignore all Poisson points at sites
$x \not \in [-L,L]^2$ without affecting the values of $(M_t(x):|x| \leq L_0)$. Hence
the final values $(m(N,x, \mathcal{Z}): |x| \leq L_0)$ will agree, on $\mathcal{G}$, with
$(M_t(x): |x| \leq L_0) $.

\vspace{.05in}
\noindent
\textbf{Claim 3.}  There are maps $K_i:\Omega_0 \to \mathcal{P}(V)$ (where $\mathcal{P}(V)$ is the power set of $V$)
so that
\begin{equation} \label{claim3.1}
A_i = \bigcup_{\omega \in A_i} [\omega]_{K_i(\omega)} \quad \mbox{for $i=1,2$,}
\end{equation}
and
\begin{equation} \label{claim3.2}
\mbox{$K_1(\omega)  \cap K_2(\omega) = \emptyset$ for all $\omega \in S^V$.}
\end{equation}
There is a natural notion of paths on $\Omega_0$. For a fixed $\omega \in \Omega_0$, and for $0 \leq k < k' \leq N$ and $|x|,|x'| \leq L$,
 we write $(k,x) \leadsto (k',x')$ if there exist $(z_j:j=0,\ldots,n) \in \Z^d \cap [-L,L]^d$ satisfying
\begin{eqnarray*}
& (1) & z_0 = x, \, z_n = x', \, n = k' - k, \\
&(2) & |x_j-x_{j-1}| \in \{0,1\} \; \mbox{for $j=1,\ldots,n$,} \\
& (3) & \mbox{if $x_j -x_{j-1} = \pm e_i$ then $\omega_{k+j,x_j}=\pm e_i$ for all $j=1,\ldots,n$,}\\
& (4) & \mbox{if $x_j -x_{j-1} = 0$ then $\omega_{k+j,x_j} \in \{0,p,q\}$ for all $j=1,\ldots,n$.}
\end{eqnarray*}
In words, the values of $\omega_{k+j,z_j}$ along $j=0,\ldots,n$ guarantee that a particle starting at $(k,x)$ would end up at $(k',x')$.
Sites $(k,x)$ can have a (unique) path to at most one site $(N,y)$ (some will exit the region $[-L,L]^d$).
Any two such paths will coalesce at the first site in common.
Now define, for $\omega \in \Omega_0$,
\[
K_i(\omega) = \left\{ (k,x): (k,x) \leadsto (N,y) \; \mbox{for some $y \in J_i$}\right\}.
\]
Then $K_i(\omega)$ is the union of the coalescing paths that lead to a site in $J_i$.
Claim (\ref{claim3.2}) is immediate since a path cannot lead both to $J_1$ and $J_2$.
It remains to check that if $\omega \in A_i$ then $[\omega]_{K_i(\omega)} \subseteq A_i$.
Observe that  $(m(k,x,\omega):x \in J_i)$ can be calculated using only on the values
$(\omega_{k,x}:(k,x) \in K_i(\omega))$, indeed $m(N,y,\omega)$ is found by tracking the mass of the initial condition and
depositions along paths that lead to $(N,y)$. For $\omega' \in [\omega]_{K_i(\omega)}$ we have
$K_i(\omega) \subseteq K_i(\omega')$ by definition of a cylinder set, with a strict inclusion possible
if other sites are now connected to $((N,y):y \in J_i)$. But these extra sites can only lead to larger values
of mass ending up in $J_i$, that is $m(N,y,\omega') \geq m(N,y,\omega)$ for all $y \in J_i$. Since $f_i$ are
non-decreasing we must have
\[
f_i \left(m(N,y,\omega'): y \in J_i \right) \geq f_i \left(m(N,y,\omega): y \in J_i \right) \geq a_i
\]
and hence $\{\mathcal{Z} = \omega'\} \cap \mathcal{G} \subseteq \hat{A}_i \cap \mathcal{G}$ as desired.

\vspace{.05in}
\noindent
\textbf{Proof of negative association.}
We proceed as in Lemma 2.4 of \cite{vandenberg+K}. Using claims 2 and 3
\begin{eqnarray*}
\hat{A}_1 \cap \hat{A}_2 \cap \mathcal{G} & \subseteq &  \mathcal{Z} \in A_1 \cap A_2 \\
& = & \mathcal{Z} \in \bigcup_{\omega \in A_1} [\omega]_{K_1(\omega)} \cap \bigcup_{\omega \in A_2}
 [\omega]_{K_2(\omega)} \\
& \subseteq &  \mathcal{Z} \in A_1 \Box A_2,
\end{eqnarray*}
The last line following from (\ref{claim3.2}).
Hence, applying BKR to the law of $\mathcal{Z}$,
\begin{eqnarray*}
P[\hat{A}_1 \cap \hat{A}_2] & \leq & P[\mathcal{Z} \in A_1 \Box A_2] + P[\mathcal{G}^c] \\
& \leq & P[\mathcal{Z} \in A_1]
P[\mathcal{Z} \in A_2] + P[\mathcal{G}^c] \\
 & \leq & P[\hat{A}_1] P[\hat{A}_2] + 3 P[\mathcal{G}^c] .
\end{eqnarray*}
Now claim 1 completes the proof.

\vspace{.05in}
\noindent
\textbf{Remark.} One might start by proving negative association for a purely discrete
time process on a  finite lattice, defined on a finite product structure, and this avoids the need for the
good set $\mathcal{G}$. One then might show that the such discrete time processes
converge to the desired continuous limit, noting the conclusion carries over naturally under
convergence in distribution.
The trick of embedding a discrete time structure inside a graphical structure for the
desired infinite lattice continuous process, and the good set $\mathcal{G}$, is designed specifically
to avoid a lengthy weak convergence argument.
\section{Asymptotics}  \label{s3}
Examining the upper bound (\ref{upperbound}) and lower bound (\ref{lowerbound}) on $q_c(p)$ for large and small
$p,q$ one finds the asymptotic inequalities
(that is in terms of limiting quotients)
\begin{eqnarray*}
& p-2\sqrt{p} \preceq  q_c(p) \preceq  p-1 & \mbox{as $p \to \infty$,} \\
& \frac{p^2}{4} \preceq  q_c(p) \preceq  \frac{p}{4} & \mbox{as $p \to 0$.}
\end{eqnarray*}
In this section we show that more accurate, dimension dependent, asymptotic information on the
critical curve is often possible.
\subsection{Large $p,q$}   \label{s3.1}
For large $q,p$ there are many evaporations and depositions between each random walk step.
The correct approximating model is one with continuous masses in $[0,\infty)$ at discrete lattice
sites $x \in \Z^d$. We first study this continuous mass model.
\subsubsection{Continuous mass model.} \label{s3.1.1}
Masses $(X_t(x): x \in \Z^d, t \geq 0)$, indexed over sites
in $\Z^d$, evolve as follows. Between jumps the masses $X_t(x)$ evolve according to the
stochastic differential equations
\begin{equation} \label{masssde}
dX(x) = - \beta \, dt + \sqrt{2}\,  dB(x) + dL(x),
\end{equation}
driven by independent Brownian motions $(B_t(x): x \in \Z^d, t \geq 0)$.
The parameter $\beta$ will reflect the excess of the evaporation rate $q$ over the deposition rate $p$
in the original model.
The term $L(x)$ is the local time of $X(x)$ at zero, so that $dL(x)$ charges only the set
$\{t: X_t(x)=0\}$ and ensures that  $X_t(x) \geq 0$.
Additionally the mass at each site $x$ has independent rate one Poisson driven jumps, where
at a jump time $t$ a nearest neighbour
$y$ is chosen at random and the masses change to
\[
X_t(x) = 0, \quad X_t(y) = X_{t-}(x) + X_{t-}(y),
\]
that is the mass at $x$ jumps onto $y$ and coalesces.
We now argue that the moment method again establishes the existence of a phase transition
with a critical parameter $\beta_c(d)$. Since we use it only as a guide to a later calculation for
the discrete mass model, we do not give a careful construction of the model or
prove the required monotonicity and dependence properties (although
one way to establish them is to take limits in the corresponding properties for discrete mass model).
Define moments for the system, with zero initial condition, by
$m_k(t) = E[X^k_t(0)]$ and set the critical value as
\[
\beta_c(d) = \inf\{ \beta \geq 0: m_1(\infty) = \infty \}.
 \]
We aim to show $\beta_c(d) \in (0,\infty)$.
 The second moment equation works smoothly. Ito calculus and negative correlation show that
\[
 \frac{dm_2(t)}{dt} =  2 - 2 \beta m_1(t)  + 2 E[X_t(0) X_t(1)] \leq 2 - 2 \beta m_1 + 2 m_1^2.
 \]
The argument from section \ref{s2.2} shows that $m_1$ stays bounded if $\beta > 2$ so that
$\beta_c(d) \leq 2$.
For the first moment one has $ \frac{d}{dt} m_1(t) = - \beta +  \frac{d}{dt} E[L_t(0)]$.
A direct estimate on the local time term does not seem useful, but we
may examine this equation at discrete times $t=1,2,\ldots$ to find
\[
m_1(n+1) - m_1(n) = - \beta + E[L_n(0)] - E[L_{n-1}(0)].
\]
The idea is to show that $E[L_t(0)] - E[L_{t-1}(0)]$ is bounded below for all $t \geq 0$ by a non-zero constant.
Then the first moment equation shows at least linear growth of the first moment if $\beta$ is small enough
and hence that $\beta_c(d) > 0$.
In the interval $[t-1,t]$ there is probability $e^{-2}$ of the event $\Omega$ that there is a single jump time
$\tau \in [t-1,t-\frac12]$ of mass away from
the origin, no further jumps away from the origin during $[t-\frac12,t]$ and
 and that no neighbouring sites jump onto the origin during $[t-1,1]$. Conditional on $\Omega$,
the local time of the mass at the origin evolves, during $[\tau,t]$, like the local time $\hat{L}$
for a (rate $2$) reflected Brownian motion $\hat{X}$ with drift $-\beta$ starting at zero, that is
\[
\hat{X}_s = W_s - \beta s + \hat{L}_s, \quad \mbox{$s \geq 0$.}
\]
Then $E[\hat{L}_s] \geq E[\hat{X}_s]$. Moroever
a comparison argument shows that $E[\hat{X}_s]$ is decreasing in $\beta$, and when
$\beta=0$ the variable $\hat{X}_s$ is just modulus of a normal $N(0,s)$ variable.
This shows there exists $c_1>0$, independent of
$\beta$, so that $ E[\hat{L}_s] \geq c_1$ for all $s \in [\frac12,1]$, establishing the desired
lower bound.
\subsubsection{Proof of the asymptotics (\ref{crude1}).} \label{s3.1.2}
We return to the discrete mass model.
Consider a single isolated site undergoing only evaporation at rate $p$ and
monomer deposition at rate $q=p-\beta p^{1/2}$. Then the rescaled mass $X_t = p^{-1/2} M_t$ is a Markov
chain satisfying, on $\{M_t \geq 1\}$,
\begin{eqnarray*}
E[ X_{t+\Delta} - X_t | \sigma(X_s:s \leq t) ] &=&  \beta \Delta + O(\Delta^2), \\
E[ |X_{t+\Delta} - X_t|^2 | \sigma(X_s:s \leq t) ] &=&  (2- \beta p^{-1/2}) \Delta + O(\Delta^2).
\end{eqnarray*}
Using these it is straightforward to establish a diffusion approximation:
the rescaled mass $X_t = p^{-1/2} M_{t}$ converge in distribution, as $p \to \infty$, to
the diffusion in (\ref{masssde}).
Extending this to the $\Z^d$ case one expects that, when $q=p-\beta p^{1/2}$,
the rescaled masses $((p^{-1/2}M_t(x): x \in \Z^d):t \geq 0)$ converge in distribution
to the continuous mass model described above.
This convergence together with the existence of the critical
value $\beta_c(d)$, suggest the asymptotics in  (\ref{largepq}).
The convergence could be set up, for example, in the space of continuous
paths taking values in some set of
measures on $\R$ (simply the empirical measure of the particle positions and masses)
with a (suitable) weak topology. However convergence for all finite time intervals $[0,T]$ would not be sufficient
to read off the asymptotic. One would like the convergence of the associated stationary
distributions (and note the convergence to equilibrium is not expected to be exponential).
We leave these ideas for later, but we quickly show that
the argument for the continuum case applies approximately to the (suitably scaled)
lattice case, although this inevitably leads to the cruder upper and lower constant
asymptotics that do not reflect the continuum critical value.
Define scaled moments by $\bar{m}_k(t) = p^{-k/2} m_k(t)$.
Choose $q=p-\beta p^{1/2}$.
In these scaled variables the the second moment inequality (\ref{qe}) becomes
\[
\bar{m}_1^2(t) - \beta \bar{m}_1(t) + (1- p^{-1/2} \beta) \geq 0 \quad \mbox{for all $t \geq 0$.}
\]
Then if $\beta >2$ and $p$ is large enough we see
that $\bar{m}_1(t)$, and hence also $m_1(t)$, remains bounded. This shows that
$\limsup_{p \to \infty} p^{-1/2}(p-q_c(p)) \leq 2$.
The scaled first moment equation (\ref{m1}) becomes
\begin{equation} \label{mm1}
\frac{d\bar{m}}{dt} = - \beta + p^{1/2} (1-s(t)).
\end{equation}
The argument for the continuum model suggests bounding $1-s(t) = P[M_t(0) =0]$ from below.
For $t \geq 1$, there is a probability $e^{-2}$ of the event $\Omega$ that
there is single jump time $\tau \in [t-1, t-\frac12]$
where the mass at the origin jumps to a neighbour, no further jumping of mass from the origin
during $[t-\frac12,t]$, and no mass jumps into the origin during $[t-1,1]$.
Conditional on $\Omega$, the mass at the origin during the interval $[\tau,t]$ follows a reflected
random walk, increasing by one at rate $q$ and decreasing by one at rate $p$,
and starting at zero.  For a bound from below
of $P[M_t(0)=0| \Omega]$ we may compare, using a simple coupling, with a reflected symmetric random
walk $\hat{M}$ at rate $p$ (that is we may, as in the continuous mass case, take $\beta = 0$).
But then exact formulae for simple symmetric random walks shows that there exists $c_1^{'}>0$
so that $P[\hat{M}_s=0] \geq c_1^{'} p^{-1/2}$ for all $ p \geq 1$ and $s \in [\frac12,1]$.
Using this estimate in (\ref{mm1}) shows that when $\beta$ is small enough
the first moment grows at a strictly positive rate, indeed that $p^{-1/2}(q_c(p) -p) \geq c_1^{'}e^{-2}$,
for all $ p \geq 1$,  completing the crude asymptotics.
\subsection{Small $p,q$ in $d=1$} \label{s3.2}
In $d=1$, for small $p,q$ the system is well approximated, on large
space-time scales, by the analogous continuous space system of massive
coalescing Brownian motions on $\R$ with evaporation and Poisson immigration.
We study this model first.
\subsubsection{Continuous space model} \label{s3.2.1}
%
Standard independent Brownian motions on $\R$ instantly coalesce upon meeting. Particles have masses
with values in $\N$, and masses add at coalescence.
Also there is deposition of monomers at a Poisson rate $q$, and mass evaporation of each particle, reducing
the mass by one, at rate $p$.
Let $(X^i_t,M^i_t: i \in I_t)$ list the positions and masses at time $t$.
We use $Q^{R,N}_{p,q}$ to denote the distribution for this model considered with zero initial conditions.
Under diffusive scaling, namely $\overline{X}^i_t = c^{-1} X^i_{c^2t}$ and $\overline{M}^i_t = M^i_{c^2t}$ for $c>0$,
we obtain a  new system of coalescing particles with distribution $Q^{R,N}_{c^2 p, c^3 q}$.
We will argue below, by similar arguments to the lattice case,
that the one parameter system, under $Q^{R,N}_{1,\alpha}$, has a single critical $\alpha_c \in (0,\infty)$.
Then the scaling shows the continuum system under $Q^{R,N}_{p,q}$
has the critical value of  $q^{R,N}_c (p) = \alpha_c p^{3/2}$.
Consider a system with distribution $Q^{R,N}_{1,\alpha}$. Since our argument for a critical value of $\alpha$
will only be used as a guide for a corresponding lattice argument,
we do not give a careful construction of the model or proof of the required monotonicity and dependence properties.
We may define the occupation density and moments by fixing a test function $\phi:\R \to [0,\infty)$ with
$\int \phi =1$ and setting
\[
s(t) = E  \sum_{i} \phi(X^i_t), \quad m_k(t) = E  \sum_i |M^i_t|^k \phi(X^i_t).
\]
Note by translation invariance that any such test function $\phi$ defines the same quantity.
Monotonicity still implies that $s(t)$ and $m_k(t)$ (for $k >0$) are non-decreasing in $t$ and $\alpha$
and we set
\[
\alpha_c = \inf\{ \alpha: m_1(\infty) =\infty \}.
\]
The first moment equation for $m_1(t)$ is
\[
\frac{dm_1}{dt}(t) = \alpha - s(t).
\]
We may bound $s(t)$ from above by $\bar{s}(t)$ the occupation density for the system with
zero evaporation, that is under the distribution $Q^{R,N}_{0,\alpha}$.
For this simpler system the scaling, as above, shows that the stationary state has a scaling property.
Indeed the stationary occupation distribution under $Q^{R,N}_{0,\alpha}$
is a scaled copy, by a factor of $\alpha^{1/3}$, of the stationary distribution
under $Q^{R,N}_{0,1}$.  Thus $\bar{s}(\infty) =  c_2 \alpha^{1/3}$
for some $c_2 \in (0,\infty)$ (in fact $c_2$ is known - see the appendix).
Substituting in $s(\infty) \leq c_2 \alpha^{1/3}$ into the first moment equation we find that
$\alpha_c \leq c_2^{3/2}$.
To find a lower bound on $\alpha_c$ we consider the second moment $m_2(t)$.
Again it seems difficult to use a differential inequality (a problem familiar from
the density decay estimates for coalescing Brownian
motions in Bramson and Griffeath \cite{Bramson+G}).
Instead we watch $m_2(t)$ along integer times $t = 1,2,\ldots$ and bound
$m_2(n+1) - m_2(n)$ from above. The argument is a little fiddly as
each of the dynamics (coalescence, evaporation and deposition) contribute, so
we start with a heuristic overview. The contribution to
$m_2(n+1) - m_2(n)$ due to coalescence will be bounded by $C m_1(n)^2$ by using
negative correlation for the mass distribution at time $n$. The contribution due to
the rate one evaporation will be bounded by $- 2 m_1(n) + s(\infty) \leq
-2m_1(n) + c_2 \alpha^{1/3}$. Finally the deposition of monomers will contribute
$2 C \alpha m_1(n) + C \alpha^2$. Monotonicity of $n \to m_2(n)$ leads to
\begin{equation} \label{diffeqn}
0 \leq m_2(n+1) - m_2(n) \leq C m_1(n)^2  + 2C \alpha m_1(n)  - 2 m_1(n) + C \alpha^2 + c_2 \alpha^{1/3}.
\end{equation}
For small $\alpha$ the quadratic inequality prevents $m_1(n)$ from entering the
interval $(\frac{1}{2C}, \frac{3}{2C})$. But
the first moment equation shows $m_1(n+1) - m_1(n) \leq \alpha$. This prevents
$m_1(n)$ jumping over this forbidden interval and forces it to remain bounded
for such small $\alpha$. Hence the critical value $\alpha_c$ is strictly positive.
Now we expand the above argument with a more careful construction of the process between times
$n$ and $n+1$. The construction exploits the fact that systems of coalescing particles can be constructed
particle by particle, with each new particle following an independent motion until it first hits
a previously constructed path, and then following the path of the particle it hits.
List the (random) positions and masses of the particles at time $n$
as $(z_i,m_i:i \in \N)$. We will give a pathwise construction of the system at time $n+1$ in three steps.

\vspace{.05in}
\noindent
\textbf{Step 1: Run CBMs.}
Run coalescing Brownian paths (independent from the evolution up to time
$n$) starting from $(z_i:i \in \N)$ over the time interval $[0,1]$.

\vspace{.05in}
\noindent
\textbf{Step 2: Add deposition.} Construct (independently) the rate $\alpha$ Poisson points
$(y_i,t_i: i \in \N)$ in $\R \times [0,1]$ of the monomer depositions.
Then construct (independently) coalescing Brownian paths started from these points that also
coalesce with any paths in step 1.
Now attach masses to the particles by tracking the mass through each coalescence. This
yields particle positions and masses $(X^i_t, \hat{M}^i_t: i \in I_t)$ for $t \in [0,1]$.

\vspace{.05in}
\noindent
\textbf{Step 3: Add evaporation.} Add (independently) Poisson marks onto the paths
of these coalescing Brownian motions to indicate the times of rate one evaporations.
Define the adjusted masses $(M^i_t : i \in I_t)$ along the paths due to these evaporations,
(so that some $M^i_t$ may take value $0$).\\
\\
After all three steps the collection $ (X^i_t,M^i_t: i \in I_t, M^i_t \neq 0)$
is then a realization of the true system at time $n+t$.
We now estimate the changes in the second moment through each of these
steps in the construction. Consider the change in the second moment after steps $1$ and $2$, namely
\[
\Delta_1 : = \sum_{i \in I_1} \phi(X^i_1) |\hat{M}^i_1|^2 - \sum_i \phi(z_i) m_i^2.
\]
We rewrite $\Delta_1$ in a different way.
For each particle $(z_i,m_i)$ we let $Z^i$ denote the final position of this
mass at time $t=1$ after step 1. Thus $Z^i$ will agree with one of the positions $(X^i_1: i \in I_1)$ but due to
coalescence we may have $Z^i = Z^j$ for many $i \neq j$. Similarly for each position $(y_i,t_i)$ we let
$Y^i$ denote final position of mass one monomer that was deposed at $(y_i,t_i)$. Then
\begin{eqnarray*}
 \sum_{i \in I} \phi(X^i_1) |\hat{M}^i_1|^2
& =  & \sum_{i \in I_1} \phi(X^i_1) \left(\sum_j m_j \chi(Z^j = X^i_1) + \sum_k \chi(Y^k = X^i_1) \right)^2 \\
& = & \sum_{j} \phi(Z^j)  m_j^2 +
\sum_{j \neq j'}  \phi(Z^j) \chi(Z^j=\Z^{j'}) m_j m_{j'} \\
& & +  \sum_{k,k'} \phi(Y^k) \chi(Y^k=Y^{k'})
+  2 \sum_{j,k} m_j \phi(Z^j) \chi(Z^j=Y^k) \\
& =:& \Delta_{1,1} + \Delta_{1,2} + \Delta_{1,3} + \Delta_{1,4}.
\end{eqnarray*}
Each position $Z_j$ is, conditionally on $\sigma((x_i,m_i:i \in \N))$, a $N(z_j,1)$ variable.
Hence, writing $(P_t)$ for the Brownian semigroup,
\[
E[ \Delta_{1,1}] = E \sum_{j \in Z}  P_1 \phi(z_j) m_j^2  = m_2(n)
\]
by translation invariance.
Each pair $(Z_j,Z_k)$, conditionally on $\sigma((x_i,m_i:i \in \N))$ and when $j \neq k$,
are the positions of a coalescing pair of Brownian motions
at time $t=1$ started at $(z_j,z_k)$. We choose $\phi(x) = \frac12 \chi(|x| \leq 1)$ and use
the following crude estimate.  Suppose $(B^x,B^y)$ denotes a pair of coalescing Brownian motions
started at $(x,s_1)$and $(y,s_2)$. Then there exists $c_3$ so that for all $0 \leq s_1 \leq s_2 \leq 1$ and $x,y \in \R$
\begin{equation} \label{quick2}
P \left[ |B^x_1| \leq 1 \mbox{ and $B^x$, $B^y$ coalesce before time $1$} \right]
\leq c_3 (|x|^{-2} \wedge 1) (|x-y|^{-2}  \wedge 1).
\end{equation}
See the appendix for a short derivation. Using this we may bound the expectation
\begin{eqnarray*}
E[ \Delta_{1,2}] & \leq &  c_3 E \sum_{j \neq j'} m_j m_{j'}
(|z_j|^{-2} \wedge 1) (|z_j-z_{j'}|^{-2} \wedge 1)  \\
& \leq & c_3 m^2_1(n) \int_{R^2} (|z|^{-2} \wedge 1) (|z-z'|^{-2} \wedge 1)  dz \, dz' \\
& = & 16 c_3 m^2_1(n).
\end{eqnarray*}
Here we have used negative correlation for the second inequality.
Next, using (\ref{quick2}) for each pair of deposed particles (arriving at Poisson $\alpha$ rate),
\[
E[ \Delta_{1,3}] \leq  c_3 \alpha^2 \int \! dy \int \! dy' \int^1_0 ds_1 \int^1_0 ds_2
(|y|^{-2} \wedge 1)(|y-y'|^{-2} \wedge 1) = 16  c_3 \alpha^2,
\]
and
\[
E[ \Delta_{1,4}] \leq 2 c_3  \alpha \int \! dy  \int^1_0 \! ds \int
 E \sum_j m_j  (|z|^{-2} \wedge 1) (|z-y|^{-2} \wedge 1) = 32 c_3 \alpha m_1(n).
\]
Next we look at the change due to evaporation defined by
\[
\Delta_2 = \sum_{i \in I_1} \phi(X^i_1) |M^i_1|^2 - \sum_{i \in I_1} \phi(X^i_1) |\hat{M}^i_1|^2.
\]
The decrease in mass due to evaporation occurs at a Poisson rate along the paths. When
a mass $m$ at $x$ suffers an evaporation it produces a jump of size
$ \phi(x) (1-2 m)$. Hence
\begin{eqnarray*}
E[ \Delta_2 ]  &=&  E \int^1_0 \sum_{i \in I^1_t}  \phi(X^i_t)
\chi(M^i_{t-} >0) dt
-  2 E \int^1_0 \sum_{i \in I^1_t}  \phi(X^i_t) M^i_{t-} dt \\
& = & \int^1_0 s(t) \, dt  - 2  \int^1_0 m_1(n+t) \, dt \\
& \leq & c_2 \alpha^{1/3} - 2 m_1(n).
\end{eqnarray*}
Combining the estimates on $\Delta_1$ and $\Delta_2$ justifies the difference inequality (\ref{diffeqn}) given at the start
with the choice $C = 16 c_3$.
\subsubsection{Proof of the asymptotics (\ref{crude2}).} \label{s3.2.2}
We return to the lattice model and choose $q= \alpha p^{3/2}$ for some $\alpha>0$.
List at time $t$ the positions and corresponding masses of all particles as
$(X^i_t,M^i_t: i \in I_t)$. Then the rescaled process
\[
\left((p^{1/2} X^i_{p^{-1} t}, M_{p^{-1}t}: i \in I_{p^{-1}t}): t \geq 0 \right)
\]
should converges in distribution as $p \downarrow 0$ to a system of coalescing massive Brownian motions
with distribution $Q^{R,N}_{1,\alpha}$ as described above. This convergence, together with the critical
value $\alpha_c$, suggest the asymptotics in (\ref{smallpq1}). However, here we simply aim for crude
asymptotics by mimicking the continuous argument for the lattice model.
Guided by the continuous space argument we may bound the occupancy
$s(t)$  by the stationary occupancy probability $\overline{s}(\infty)$  for the zero evaporation model.
For $d=1$ models without evaporation a (time reversal) Markov duality tool is applicable and
one finds $\overline{s}(\infty) = 1- h_q(1)$ where
\begin{equation} \label{duality}
h_q(x) = E  \exp \left(-q \int^{\tau}_0 X^x_s ds\right), \quad \mbox{with
$\tau = \inf\{t: X^x_t = 0\}$,}
\end{equation}
and $(X^x_t:t\geq 0)$ is a rate $2$ continuous time simple random walk on $\N$ started at $x$.
Details of the duality argument leading to (\ref{duality}) are in the appendix \ref{s4.3}.
Note that $h_q \in [0,1]$ solves
\[
\Delta h_q(x) =   q x h_q(x) \quad \mbox{for $x \in \N_+ = \{1,2,\ldots\}$,}
\]
where $\Delta$ is the discrete Laplacian on $\N$,
with boundary conditions $h_q(0) =1$ and $h_q(x) \to  0$ as $x \to \infty$.
To see the scaling of $h_q(1)$ in $q$ we define $A_q(x) = h_q(xq^{-1/3})$, so that
$A_q$ solves
\[
\Delta A_q(x) =  x A_q (x) \quad \mbox{for $x \in q^{1/3} \N_+$,}
\]
where $\Delta$ is the discrete Laplacian on $q^{1/3} \N$,
with boundary conditions $A_q(0) =1$ and $A_q(x) \to  0$ as $x \to \infty$.
Then one can check that $\sup_{x}|A_q(x) - A_0(x)| \to 0$ as $q \to 0$, where
$A_0:[0,\infty) \to [0,1]$ is the corresponding  $C^2([0,\infty))$ Airy function solving $A^{''}_0(x) = x A_0(x)$.
Moreover one can check (details are in the appendix \ref{s4.4}) that the discrete derivative
\begin{equation} \label{Airy}
q^{-1/3} \left(1-h_q(1)\right) =  q^{-1/3}\left(A_q(0) - A_q(q^{1/3})\right) \to  - A^{'}_0(0) >0,
\end{equation}
Substituting $s(t) \leq \overline{s}(\infty) = 1- h_q(1)$ into the first moment equation
$dm_1/dt = q - ps(t)$ and using the above asymptotics leads to
\[
\limsup_{p \to 0} q_c(p) p^{-3/2} \leq (-A_0^{'}(0))^{3/2}.
\]
For the second moment equation we again mimic the argument for the continuous space Brownian motion case.
We consider the second moments $m_2(t)$  at the discrete set of times $t_n = n p^{-1}$.
The scaling above suggests that $m_2(tp^{-1})$ is $O(p^{1/2})$, and so we aim for a difference inequality of the form:
there exists $C< \infty, p_0 >0$ so that for $p \leq p_0, \alpha \leq 1$
\begin{equation} \label{modified}
0 \leq  m_2(t_{n+1}) - m_2(t_n)  \leq C p^{-1/2} m_1(t_n)^2  + 2C \alpha m_1(t_n)  - 2 m_1(t_n) + C (\alpha^2 +  \alpha^{1/3}) p^{1/2}.
\end{equation}
This implies for small enough $\alpha$ that the values of $p^{-1/2} m_1(t_n)$ cannot lie in $(\frac{1}{2C}, \frac{3}{2C})$.
However the first moment equation implies that $p^{-1/2} m_1(t_{n+1}) - p^{-1/2} m_1(t_n) \leq \alpha$. So for small enough $\alpha$ the first moment is bounded. This gives
\[
\liminf_{p \to 0} q_c(p) p^{-3/2} >0.
\]
To establish (\ref{modified}) we way repeat the three step construction from section \ref{s3.2.1} to construct the process
between times $t_n$ and $t_{n+1}$ making the natural changes, namely:
replace coalescing Brownian motions on $\R$ by rate one coalescing
random walks on $\Z$; replace rate $\alpha$ deposition on $[0,1] \times \R$ by rate $q$ deposition on $[0,1] \times \Z$;
and replace rate $1$ evaporation along the paths by rate $p$ evaporation.
We may rewite the moments as $m_k(t) = E \sum_{x} \phi(x) M_t^k(x)$ for any test function
$\phi$ satisfying $ \sum_x \phi(x) = 1$. We choose $\phi(x) = N_p^{-1} \chi(|x| \leq p^{-1/2})$
where $N_p \approx 2 p^{-1/2}$ is chosen so that the constraint holds.
Decompose $m_2(t_{n+1}) - m_2(t_n)$ into terms $\Delta_1, \Delta_2$ as before, where we work however
over a time interval of length $p^{-1}$. We use an analogue of (\ref{quick2}),
which is derived in the appendix \ref{s4.2}.
Let $(X^x,X^y)$ be a pair of coalescing rate one simple random walks
started at $(x,s_1)$and $(y,s_2)$. Then there exists $c_3^{'}$ so that for all $0 \leq s_1 \leq s_2 \leq p^{-1}$,
$x,y \in \Z$ and $p \leq 1$
\begin{equation} \label{quick3}
P \left[ |X^x_{p^{-1}}| \leq p^{-1/2} \mbox{and $X^x$, $X^y$ coalesce by time $p^{-1}$} \right]
\leq c_3^{'} (p^{-1} |x|^{-2} \wedge 1) (p^{-1} |x-y|^{-2}  \wedge 1).
\end{equation}
Following the steps from  section \ref{s3.2.1}, we find
$E[\Delta_{1,1}] = m_2^2(t_m)$ and
\begin{eqnarray*}
E[\Delta_{1,2}] & \leq & c_3^{'} N_p^{-1} E \sum_{j \neq j'}
m_j m_{j'} (p^{-1} |x|^{-2} \wedge 1) (p^{-1} |x-y|^{-2}  \wedge 1) \\
& \leq & c_3^{'} N_p^{-1} m_1^2(t_n) \sum_{z,z'}
 (p^{-1} |z|^{-2} \wedge 1) (p^{-1} |z-z'|^{-2}  \wedge 1) \\
& \leq & \frac13 c_3^{'} c_4^2 m_1^2(t_n) p^{-1/2}.
\end{eqnarray*}
The second inequality uses the negative correlation in (\ref{pairwisenc}), and the third inequality uses
$N_p \leq  p^{-1/2}$ and $\sum_z (p^{-1} |z|^{-2} \wedge 1) \leq c_4 p^{-1/2}$ for all $p \leq 1$.
In a similar way we find
\[
E[\Delta_{1,3}] \leq  \frac13 c_3^{'} c_4^2  \alpha^2 p^{1/2} \quad \mbox{ and} \quad
E[\Delta_{1,4}] \leq  \frac23 c_3^{'} c_4^2 \alpha  m_1(t_m).
\]
Using (\ref{Airy}) we may choose $p_0 \leq 1$ small enough that
\[
s(\infty) \leq (-2 A^{'}_0(0)) q^{1/3} = (-2 A^{'}_0(0)) \alpha^{1/3} p^{1/2}.
\]
Then $E[\Delta_2] \leq s(\infty) - 2 m_1(t_n) \leq (-2 A^{'}_0(0))  \alpha^{1/3} p^{1/2} - 2 m_1(t_n)$.
Together these estimates imply (\ref{modified}) with
$C = (-2 A^{'}_0(0))  \vee   \frac13 c_3^{'} c_4^2$, completing the proof.
\subsection{Remarks on small $p,q$ in $d \geq 3$}   \label{s3.3}
There is no continuous space model approximation in $d>1$, and the intuition is somewhat different.
We give a purely heuristic argument suggesting that $q_c(p) \approx p^2$ as $p,q \downarrow 0$.
The argument exploits the small occupation densities and the transience of random walks
in $d \geq 3$. It follows the intuitive argument given in \cite{vandenberg+K}
for the modified rate equations that gives the decay of the occupation density in coalescing random walks.
In  \cite{vandenberg+K} they exploit large times to obtain low densities, where we exploit small
$p,q$ values.
We use the informal notation $z \to z' $ to denote the event that
the particle at $z$ moves to position $z'$.
Choose $0 << r << t$ so that $p^{-2/d} <<  r << p^{-1}$. Choose $q = O(p^2)$. The first moment equation
shows that the occupation density satisfies $s(t) \leq q/p = O(p)$. Then the following approximations seem
reasonable for small $p$:
\begin{eqnarray}
 P[M_t(0)>0, \,M_t(e) >0]
& \approx & P[ \cup_{x \neq y} \{M_{t-r}(x)>0, \, M_{t-r}(y) >0, \, x \to 0, y \to e\}] \nonumber  \\
& \approx & \sum_{x \neq y} P[ M_{t-r}(x)>0, \, M_{t-r}(y) >0, \, x \to 0, y \to e] \nonumber  \\
& \approx & \sum_{x \neq y} p_d  p_r(x,0) p_r(y,e) P[ M_{t-r}(x)>0, \, M_{t-r}(y) >0] \nonumber  \\
& \approx & \sum_{x, y} p_d  p_r(x,0) p_r(y,0) P[ M_{t-r}(x)>0] P[M_{t-r}(y) >0] \nonumber  \\
& = & p_d s^2(t-r) \nonumber  \\
& \approx &  p_d s^2(t). \label{mmf}
\end{eqnarray}
We expect all expressions in these approximations to be $O(p^2)$.
The first approximation throws away the event that the particles at $0,e$ were
deposited between $t-r$ and $t$. The expected mass of deposited particles
during $[t-r,t]$ that reach $0$ at time $t$ is of order $q r << p$, suggesting this
approximation is valid.
The second approximation uses low densities: the error should involve three occupied
sites and be of order $O(p^3)$.
The third approximations involves the probability
$P[x \to 0, y \to e \, \mbox{without coalescing}]$. Evaporation at rate $p$ can be ignored
for these two particles over this interval since $r << p^{-1}$. Considering the time reversed paths we need
to find $P[0 \to x, e \to y \, \mbox{without coalescing}]$, which by transience of the walks is well approximated,
for $|x-y| >> 1$, by $p_d p_r(0,x) p_r(0,y)$.
Negative correlation shows that the fourth approximation is an upper bound. Moreover
\begin{eqnarray*}
&& \hspace{-.3in} \sum_{x \neq y}  p_r(x,0) p_r(y,e) P[ M_{t-r}(x)>0, \, M_{t-r}(y) >0] \nonumber
\\
& = & E\left[ (\sum_x p_r(x,0) \chi(M_{t-r}(x)>0))^2 \right] - \sum_x p^2_r(x,0) P[M_{t-r}(x)>0]
\nonumber  \\
& \geq &
 \left( E\left[ \sum_x p_r(x,0) \chi(M_{t-r}(x)>0) \right] \right)^2 - s(t-r)  \sum_x p^2_r(x,0)
\nonumber \\
& = & s^2(t-r) - s(t-r) \sum_x p^2_r(x,0).
\end{eqnarray*}
The error in the fourth approximation is thus bounded by
\[
s(t-r) \sum_{x}  p^2_r(x,0) \leq C p r^{-d/2}
\]
which is $o(p^2)$ since $p^{-2/d} << r$. Similar approximations are used in \cite{vandenberg+K},
where careful error bounds are established.
We show that just the approximation (\ref{mmf}) suggests the scaling for $q_c(p)$.
Note
\begin{eqnarray*}
0 \leq \frac{ds(t)}{dt}(t) & = & q(1-s) - p P[M_t(0)=1] - P[M_t(0)>0, \, M_t(e)>0] \\
& \leq & q -  P[M_t(0)>0, \, M_t(e)>0].
\end{eqnarray*}
Using the approximation (\ref{mmf}) above suggests that
$s(\infty) \leq p_d^{-1/2} q^{1/2}$. Using this in the first moment equation then suggests
that $\limsup_{ p \to 0} q_c(p) p^{-2} \leq p_d^{-1} p^2$. This compliments the rigourous upper bound, but we
believe we can guess the correct asymptotic.
For the case of purely coalescing random walks in $d \geq 3$, with initially one
monomer per site, the occupation density drops as $p_d^{-1} t^{-1}$.  Arratia \cite{arratia} showed  that one
can rescale space by $t^{-1/d}$ to leave an asymptotically constant density of occupied sites, and the positions
of the rescaled sites converge to a Poisson point process on $\R^d$. Moroever the mass of
the rescaled particles become independent, leading to a compound Poisson limit for the empirical measure.
The aim is to mimic this for a rescaled version of the stationary distribution of our model as $p \to 0$. There is
similar underlying intuition for a Poisson limit, as follows. Run the process in its stationary distribution for an
interval of length $t$. If $t << p^{-1}$ then in bounded regions of space one expects
negligible amounts of deposition or evaporation.
Provided $t >> p^{-2/d}$ the particles will move far beyond the
typical inter-particle distance $O(p^{1/d})$. Due to transience the process is therefore well
approximated by free independent motion over this interval, and free motion
converges to a Poisson limit.
To find the approximating compound Poisson approximation we
use a more general version of approximations such as (\ref{mmf}).
Take $q= \alpha p^2$ for some $\alpha >0$. We expect
\begin{equation}
\label{mmf2}
E[F(M_{\infty}(0)) F(M_{\infty}(x))] \approx p_d(x) E[F(M_{\infty}(0))] \, E[F(M_{\infty}(x))]
\end{equation}
for all bounded $F:\N \to \R$ satisfying $F(0) = 0$,
where $p_d(x)$ is the probability of a simple random walker never hitting the origin starting at $x$.
Taking $F(m)=z^m -1$ allow us to characterize the mass via its transform, and already suggests
that the set of massed particles will also converge to a compound Poisson limit.
To find the Poisson intensity $\hat{s}_{\alpha}$ and the generating function $\hat{\phi}_{\alpha}(z)$
for the attached masses, we start by developing
$ \frac{d}{dt} E z^{M_t(0)} $ and taking $t \to \infty$ to find in the stationarity
\[
0 = q (z-1) E z^{M_{\infty} (0)} + p (z^{-1}-1) E z^{M_{\infty} (0)} \chi(M_t(0)>0)
+ 1 - 2 E z^{M_{\infty}(0)} + E z^{M_{\infty}(0)+M_{\infty}(e)}.
\]
Define
$ \phi(z) : = E[ z^{M_{\infty}(0)} | M_{\infty}(0)>0] $
so that
\[
E z^{M_{\infty} (0)} = E e^{-\theta M_{\infty} (0)} \chi(M_t(0)>0) + P[M_{\infty}(0)=0]
= s(\infty) \phi(\theta) + 1-s(\infty).
\]
Apply (\ref{mmf2}) to approximate for small $p$
\[
 1 - 2 E z^{M_{\infty}(0)} + E z^{M_{\infty}(0)+M_{\infty}(e)} \approx
p_d \left( 1- E z^{M_{\infty} (0)} \right)^2.
\]
We find, for small $p$,
\[
0 \approx q (z-1) (s(\infty) \phi(z) +1 -s(\infty))
  + p (z^{-1}-1)  s(\infty) \phi(z)
+ p_d s^2(\infty) (1-\phi(z))^2.
\]
Now we suppose, as $p \to 0$ with $q = \alpha p^2$, that
\[
p^{-1} s(\infty) \to \hat{s}_{\alpha}, \quad \phi(z) \to \hat{\phi}_{\alpha}(z)
\]
and obtain
\[
0 = \alpha (z-1) + (z^{-1}-1)  \hat{s}_{\alpha}  \hat{\phi}_{\alpha}(z)
+ p_d \hat{s}^2_{\alpha} (1-\hat{\phi}_{\alpha}(z))^2.
\]
Solving, with the correct sign, we find
\[
2 p_d \hat{s}_{\alpha} z (\hat{\phi}_{\alpha}(z) -1) =  - (1-z)
+ \sqrt{(1-z)Q(z)}
\]
where $Q(z) =  1 -  (4 p_d \alpha \hat{s}_{\alpha} +1)z + 4 p_d \alpha z^2$.
We need to discover the unknown $\hat{s}_{\alpha}$.
We now follow the steps of the mean field analysis in \cite{MKB1} (although our case is slightly simpler
in that $Q$ is quadratic, due to reducing to a one parameter model).
The tail of the variable $M$ with generator
$\hat{\phi}_{\alpha}$ can be found by inverting the generating function via the Cauchy integral
\[
P[M=m] = (2 \pi i)^{-1} \int_{C} \frac{ \hat{\phi}_{\alpha}(z)}{z^{m+1}} dz.
\]
where $C$ is a contour around the origin.
The contour can be deformed so that it goes around the singularities of
$\hat{\phi}_{\alpha}$.
We expect a growth phase where the first moment predicted by $\hat{\phi}_{\alpha}$ is infinite and
$s_{\alpha}=\alpha$. Solving for the roots of $Q(z)$ when $s_{\alpha}=\alpha$ gives roots
$z=1$ and $z=(4 p_d \alpha)^{-1}$. Let $\alpha_c = (4 p_d)^{-1}$. Then, for $\alpha > \alpha_c$,
 $\hat{\phi}_{\alpha}$ has only a single singularity due to the branch cut for the root
$z=(4 p_d \alpha)^{-1} >1 $, and the inversion then shows the tail
$P[M=m]$ drops exponentially. For $\alpha < \alpha_c$ we consider values of
$\hat{s}_{\alpha} \in [0,\alpha]$. Looking at the sum of the roots
of $Q(z)$ one finds that there must be a root with a real part less than $1$
and the inversion would then predict exponential growth of $P[M=m]$, except if there
is a double root, which would not lead to a singularity for $\hat{\phi}_{\alpha}$.
The unique value for $\hat{s}_{\alpha}$ leading to a double root is
$(4 \sqrt{p_d \alpha} -1)/4 p_d$. We summarize below our beliefs.

\vspace{.1in}
\noindent
\textbf{Conjecture.}
Under the stationary distribution with $q= \alpha p^2$,
the scaled empirical measures
$ \sum_{x} M_{\infty}(x) \delta_{p^{1/d}x} $
converge in distribution as $p \downarrow 0$,
 in the space of locally finite point measures with the vague topology,
to a compound Poisson limit $\sum_i M_i \delta_{x_i}$, where
where $(X_i)$ are Poisson with intensity  $\hat{s}_{\alpha}$ and the masses $(M_i)$ are, conditionally, independent identically
distributed with Laplace transform $\hat{\phi}_{\alpha}(\theta)$.
\section{Appendix: details for some proofs.} \label{s4}
\subsection{Proof of Proposition \ref{basics} and Lemma \ref{monotonicity}.} \label{s4.1}
We first consider the comparison theorem and  uniqueness for solutions to (\ref{sde}).
Suppose $M,\overline{M}$ are two solutions, with the same driving Poisson processes,
satisfying $M_0 \leq \overline{M}_0$ a.s. Then
\begin{eqnarray*}
&& \hspace{-.3in} d \chi(M_t(x) > \overline{M}_t(x)) \\
& = &
- \chi ( M_{t-}(x) =1,\, \overline{M}_{t-}(x) =0) dP^{(p)}_t
-  \chi(M_{t-}(x) > \overline{M}_{t-}(x)) \sum_{y \sim x} dP_t(x,y) \\
& & + \sum_{y \sim x} \left(  \chi(M_{t-}(x) +  M_{t-}(y) > \overline{M}_{t-}(x) + \overline{M}_{t-}(y))
 -  \chi(M_{t-}(x) > \overline{M}_{t-}(x)) \right) dP_t(y,x) \\
& \leq & \sum_{y \sim x}  \chi(M_{t-}(y) - \overline{M}_{t-}(y) >0) dP_t(y,x)
\end{eqnarray*}
where in the final inequality we have used $\chi(x+y>0) - \chi(x>0) \leq \chi(y >0)$.
Then, for $\theta>0$,
\[
E  \sum_x e^{-\theta |x|} \chi(M_t(x)  > \overline{M}_t(x))
\leq e^{\theta} \int^t_0  E  \sum_x e^{\theta |x|} \chi(M_s(x)  > \overline{M}_s(x))  ds
\]
and Gronwall's inequality implies that $M_t \leq \overline{M}_t$ a.s. This implies the comparison theorem
and pathwise uniqueness of solutions. As usual, uniqueness implies the Markov property.
Let $M$ be the solution started from the minimal initial condition $M_0=0$
and write $\nu_s$ for the distribution of $M_s$. The comparison theorem allows one to run a larger coupled solution
$\overline{M}$ started from a random initial condition $\overline{M}_0$ with distribution $\nu_s$
(but independent of the Poisson drivers).  Then for non-decreasing $F$ one has
\[
E\left[ F(M_t) \right] \leq E\left[ F(\overline{M}_t) \right] =
E\left[ F(M_{t+s}) \right]
\]
by the Markov property. This establishes the monotonicity in $t$ in Lemma \ref{monotonicity}.
A simpler argument establishes the monotonicity in $p$ and $q$. Let
$M$ be a solution to (\ref{sde}) and let $\overline{M}$ be a solution with the drivers
$dP^{(q)}(x)$ replaced by $dP^{(q)}(x) + dP^{(\overline{q}-q)}(x)$, where
$(P^{(\overline{q}-q)}(x): x \in \Z^d)$ is a further independent family of rate
$\overline{q}-q>0$ Poisson processes. Another Gronwall estimate as above shows, if
$M_0 = \overline{M}_0$ a.s., that $M_t \leq \overline{M}_t$ and this establishes the monotonicity in $q$.
From zero initial conditions one may check, as before, that $E[\overline{M}_t(0) - M_t(0)]$ is non-decreasing, and so
\[
0 \leq \frac{d}{dt}E[\overline{M}_t(0) - M_t(0)] = (\overline{q} - q) - p P[M_t(0)=0, \overline{M}_t(0)>0].
\]
Writing $s^{(q)}(t), s^{(\overline{q})}(t)$ for the occupancy probabilities of $M_t, \overline{M}_t$, we
find
\[
 s^{(\overline{q})}(t) - s^{(q)}(t) =  P[M_t(0)=0, \overline{M}_t(0)>0] \leq (\overline{q}-q)/p.
\]
Letting $t \to \infty$ one deduces, as in section \ref{s2.1}, that $s^{(q)}(\infty)  <q/p$ for
$s > \hat{q}_c(p)$.
Monotonicity in $p$ is similar, by constructing the evaporation driver as the sum of two
independent families $P^{(p)}(x) = P^{\overline{p}}(x) + P^{(p-\overline{p})}(x)$ of rates
$p$ and $p-\overline{p}>0$.
Existence of solutions can be established by the usual iteration argument:
set $M^{(0)}_t = M_0$ for all $t \geq 0$ and define cadlag adapted processes
by
\begin{eqnarray*}
dM^{(n)}_t(x) &=&  - \sum_{y \sim x} M^{(n-1}_{t-}(x) \, dP_t(x,y)
+ \sum_{y \sim x} M^{(n-1}_{t-}(y) \, dP_t(y,x) \\
&& \hspace{.3in}
dP^{(q)}_t(x) - \chi (M^{(n-1)}_{t-}(x) >0) \, dP^{(p)}_t(x).
\end{eqnarray*}
Note that $d(\eta,\eta') = \sum_x e^{-\theta |x|} \chi(\eta(x) \neq \eta'(x))$
is a complete metric generating the topology on $\mathcal{S}$. Using the
simple inequality $\chi(x+y \neq x'+y') - \chi(x \neq x') \leq \chi(y \neq y')$ one finds
\begin{eqnarray*}
 \sup_{s \leq t} \chi (M^{(n)}_s(x) \neq M^{(n-1)}_s(x))
& \leq & \int^t_0 \chi( M^{(n-1)}_{s-}(x) \neq M^{(n-2)}_{s-}(x)) ( dP^{(p)}_s
+ \sum_{y \sim x} dP_s(x,y) ) \\
&& + \int^t_0  \chi( M^{(n-1)}_{s-}(y) \neq M^{(n-2)}_{s-}(y)) \sum_{y \sim x} dP_s(y,x)
\end{eqnarray*}
and hence
\[
E \left[ \sup_{s \leq t} d(M^{(n)}_s,M^{(n-1)}_s) \right] \leq (2 + e^{\theta})
\int^t_0 E \left[ d(M^{(n-1)}_s,M^{(n-2)}_s) \right] ds.
\]
Now the usual argument shows that $M^{(n)}$
converges uniformly over compact time intervals to a cadlag adapted limit
 which is the solution to (\ref{sde}) started at $M_0$. An argument as for the comparison theorem shows
the solutions are non-negative. For two solutions  $M,\, M'$ with deterministic
conditions $M_0 = \eta$, $M_0^{'}=\eta'$, one has similarly
\[
E \left[ d(M_t,M^{'}_t) \right] \leq C(t,\theta) d(\eta, \eta')
\]
and this implies the Feller property for the semigroup.
The moments (\ref{moments})
are established in the same manner, that is by developing $dM^k_t(x)$
and establishing a Gronwall estimate.
For the attractiveness of the stationary distribution it is convenient to construct an extended
system that allows particles of infinite mass.  Let $\N \cup \{\infty\}$ be $\N$ with $\infty$ added as
a discrete site. At the slight risk of confusion we continue to
write the equations (\ref{sde}) for this extended system, but where we assign natural rules for
sites with the value $\infty$.
We extend the state space to
$\hat{\mathcal{S}} = (\N \cup \{\infty\})^{\Z^d}$.
We allow addition, subtraction and comparison on $\N \cup \{\infty\}$ via
\[
\infty \pm n = \infty, \quad \infty + \infty = \infty, \quad \infty - \infty = 0, \quad n < \infty.
\]
(The operation $\infty - \infty$ will only occur when a particle of mass leaves a site due to a random
walk step.)
The metric $d(\eta,\eta')$ extends to $\hat{\mathcal{S}}$ and still generates the
product topology. The proof of existence, comparison and uniqueness
of solutions carries over with only the (small) natural changes. The system is a true extension,
in that solutions to the original system continue to be solutions to the extended system, but
which never take the value $\infty$.
The comparison theorem shows there is maximal solution
$\hat{M}$ with initial condition $\hat{M}_0(x) = \infty$ for all $x$. Define
\[
A_t = \{x \in \Z^d: \hat{M}_t(x) = \infty\}.
\]
Examining the evolution of $d \chi(\hat{M}_t(x) = \infty)$ shows that
the evolution of $A_t$ is precisely that of the set of occupied sites for
instantly coalescing random walkers, started from all sites in $\Z^d$ occupied.
In particular $P(x \in A_t) $ is constant in $x$ and converges to $0$ as $t \to \infty$.
If $M$ is the minimal solution started from $M_0=0$, then comparison implies
$M_t \leq \hat{M}_t$ for all $t$. We claim that $M_t(x) = \hat{M}_t(x)$ for $x \in A^c_t$.
The proof is another simple Gronwall estimate showing
$E [ \sum_{x \in A_t^c} e^{-\theta |x|} \chi(M_t(x) < \hat{M}_t(x)) ] =0 $, by using
\[
d \chi( M_t(x) < \hat{M}_t(x) < \infty)  \leq \sum_{y \sim x}  \chi( M_{t-}(y) < \hat{M}_{t-}(y) < \infty)
dP_t(y,x).
\]
This gives all the ingredients needed to complete the proof of Proposition \ref{basics}. Starting from zero
initial conditions, the monotonicity
of the Laplace functional implies convergence in distribution to $M_{\infty} \in \hat{S}$. But the coupling
above by $\hat{M}$ and the decay of $P[x \in A_t]$ implies that $M_{\infty} \in S$ almost surely. The limit distribution
$\nu_{\infty}$ is a stationary distribution (an argument that uses the Feller property).
Any other solution can be coupled between
$M$ and $\hat{M}$ and so converges in distribution to $\nu$. In particular the stationary distribution is unique.
\subsection{Proof of the coalescent pair estimates (\ref{quick2}) and (\ref{quick3}).} \label{s4.2}
We give the argument for the discrete walk case (\ref{quick3}), and
note the Brownian case can be treated entirely similarly. We do not seek accurate bounds.
Split the problem in two by using Cauchy Schwarz, namely
\begin{eqnarray}
&& \hspace{-.6in}
\left( P [ |X^x_{p^{-1}}| \leq p^{-1/2}\mbox{ and $X^x$, $X^y$ coalesce by time $p^{-1}$}]\right)^2 \nonumber \\
& \leq &  P [ |X^x_{p^{-1}}| \leq p^{-1/2}] \,  P\left[\mbox{$X^x$, $X^y$ coalesce by time $p^{-1}$} \right]. \label{CS}
\end{eqnarray}
For $|x| \geq 2 p^{-1/2}$ and $p \leq 1$ we have
\[
P [ |X^x_{p^{-1}}| \leq p^{-1/2}]  \leq P [ |X^0_{p^{-1}}| \geq |x|/2] \leq 16 |x|^{-4} E[ |X^0_{p^{-1}}|^4]
= 16 |x|^{-4}(p^{-1} + 3 p^{-2}).
\]
Thus, for suitable $C$,  $P [ |X^x_{p^{-1}}| \leq p^{-1/2}] \leq C (p^{-2}|x|^{-4} \wedge 1)$ for all $x$ and $p \leq 1$.
Take walks $X^z$ and $X^y$ starting at $(z,s_2)$ and $(y,s_2)$. Write $Q_x$ for the law of a simple random walk
on $\Z$ started at $x$ on pathspace.
\begin{eqnarray*}
&& \hspace{-.4in} P [\mbox{$X^z$, $X^y$ coalesce by time $p^{-1}$} ] \\
& = &  Q_{z-y}[ \mbox{$X$ hits zero before $2(p^{-1}-s_2)$}] \\
& \leq &   Q_{|z-y|}[ \mbox{$X$ hits zero before $2p^{-1}$}] \\
& \leq & 2 Q_{|z-y|}[ X_{2 p^{-1}} \leq 0 ] \qquad \mbox{by the reflection principle}\\
& \leq & 2 \min\{  Q_{z-y}[ X_{2 p^{-1}} \leq 0 ],  Q_{y-z}[ X_{2 p^{-1}} \leq 0 ] \}.
\end{eqnarray*}
By conditioning at time $s_2$ we have
\begin{eqnarray*}
&& \hspace{-.4in} P [\mbox{$X^x$, $X^y$ coalesce by time $p^{-1}$} ] \\
& \leq & 2  E \min\{  Q_{X_{s_2}^x-y}[ X_{2 p^{-1}} \leq 0 ],  Q_{y-X^x_{s_2}}[ X_{2 p^{-1}} \leq 0 ] \} \\
& \leq & 2  \min\{ E Q_{X_{s_2}^x-y}[ X_{2 p^{-1}} \leq 0 ], E Q_{y-X^x_{s_2}}[ X_{2 p^{-1}} \leq 0 ] \} \\
& = & 2 Q_{|x-y|}[X_{2p^{-1} + s_2-s_1} \leq 0]\\
& \leq & C (|x-y|^{-4} p^{-2} \wedge 1)
\end{eqnarray*}
by a Markov inequality as before. Substituting these estimates into (\ref{CS}) completes the proof.
\subsection{Proof of the duality (\ref{duality}).} \label{s4.3}
The Brownian web and the dual Brownian web (see \cite{web}) is a useful construction for the study of the model in
$d=1$, at least with zero evaporation. Informally, a realization of the web on $[0,t] \times \R$ has coalescing Brownian paths
starting from every space-time point $(t,x)$. If one superimposes an independent Poisson $\alpha$ process on
$[0,t] \times \R$, yielding the space-time points $(t_i,x_i)_{i \in \N}$ of the deposed monomers, then
almost surely the web has a unique set of coalescing paths starting at each $(t_i,x_i)$. These
can be used to construct the process over $[0,t]$ under the distribution $Q^{R,N}_{0,\alpha}$. The dual web contains
an additional set of coalescing paths running backwards over the interval $[0,t]$, starting at every
point in $\{t\} \times \Q$. Moreover, almost surely, a forwards path, started from $(t_i,x_i)$, never crosses
any of the backwards dual paths. This property implies the following identity for the probability, under
$Q^{R,N}_{0,\alpha}$, of an empty interval $[x,y]$ at time $t$, where $x<y$, is given by
\[
P \left[ \mbox{$[x,y]$ is empty at time $t$} \right] = E \left[ \mbox{$[B^x_s,B^y_s] \cap (t_i,x_i)_{i \in \N} =
 \emptyset$ for all $s \in [0,t]$} \right]
\]
where $B^x_s,B^y_s$ are the positions at time $s$ of the coalescing backwards dual paths started at $B^x_t=x$ and $B^y_t=y$.
The point is that if any particle is deposited in the interval $(B^x_s,B^y_s)$ then its forwards path would have to
end up in the interval $[x,y]$ at time $t$. Taking expectations over the Poisson deposition yields
\[
P \left[ \mbox{$[x,y]$ is empty at time $t$} \right] = E  \exp\left(-\alpha \int^t_0 (B^x_s-B^y_s) ds\right).
\]
Letting $t \to \infty$ one finds the empty interval probability for the stationary distribution under $Q^{R,N}_{0,\alpha}$
is given by
\[
P \left[ \mbox{$[x,y]$ is empty at time $\infty$} \right]  =  E \exp \left(-\alpha \int^{\tau}_0 X^{y-x}_s ds\right).
\]
where $X^{y-x}$ is a rate $2$ Brownian motion started at $y-x$ and $\tau = \inf\{s: X^{y-x}_s =0\}$.
The right hand side is given by the Airy function $A_0(\alpha^{1/3}(y-x))$. Differentiating in $x$ and then letting $y \downarrow x$ gives
the stationary occupation density as  $- \alpha^{1/2} A_0^{'}(0)$. This identifies the constant $c_2$ in section \ref{s3.1}
as $c_2 = - A_0^{'}(0)$. One can see many related arguments in \cite{Pfaff}. A more detailed
description of the evolution and stationary distribution under $Q^{R,N}_{0,\alpha}$ is in preparation in a current thesis.
A continuous time random walk analogue of some of the above construction is possible,
quite close to the graphical construction of the voter model and its dual. Start with a graphical construction of
coalescing rate one simple random walks on the lattice $[0,t] \times \Z$ as follows.
Let $(t_i,y_i,\pm)_{i \in \N}$ be the points of a rate one marked Poisson process on the dual lattice
$[0,t] \times (\Z + \frac12)$. The marks $\pm$ are chosen independently for each point, with equal probability.
Passing a point $(t,n+\frac12,+)$ a forwards path at $n$ will jump to $n+1$.
Passing a point $(t,n+\frac12,-)$ a forwards path at $n+1$ will jump to $n$.
These rules define a unique right continuous
forwards paths started at any point in $[0,t] \times \Z$. Moreover paths started at different points will coalesce
instantaneously should they meet.
The same marked points can be used to define backwards paths living on the dual lattice.
When such a path approaches a point  $(t,n+\frac12,+)$ it will jump to $n-\frac12$, and
when it approaches a point  $(t,n+\frac12,-)$ it will jump to $n+\frac12$.
These rules define right continuous paths $(Y^{x+\frac12}_s: x \in \Z, s \in [0,t])$ so that
their backwards versions $s \to Y^x_{t-s}$ are coalescing  rate one random walks started at $X^{x+\frac12}_t
= x+ \frac12$.  One can then check that
forwards paths never jump over backwards paths in the following sense: if a forwards path $X$ and a backwards path
$Y$ are both defined at $s$ and $t$, then
\[
X_s < Y_s \quad \mbox{if and only if} \quad X_t < Y_t.
\]
Now the argument is as for the Brownian web.
Lay down an independent set of deposition points $(t_i,x_i)_{i \in \N}$ as a rate $q$ Poisson process on  $[0,t] \times \Z$,
and use the forwards paths constructed above to build the process on $[0,t]$ that has distribution $Q^{Z,N}_{0,q}$.
A site $x \in \Z$ is empty at time $t$ precisely if the interval $[Y^{x+\frac12}_s,Y^{x-\frac12}_s]$
contains none of the deposition points $(t_i,x_i)_{i \in \N}$ at any time $s \in [0,t]$.
Now mimic the Brownian argument to reach (\ref{duality}).
\subsection{Proof of convergence to the Airy function and (\ref{Airy}).} \label{s4.4}
Note that $A_q$ and $A_0$ are convex, non-increasing and take values in $[0,1]$.
Let $DA_q(x) = q^{-1/3}(A_q(x+q^{1/3})-A_q(x))$ be the discrete derivative, so that
 the discrete Laplacian is $\Delta A_q(x) = q^{-1/3}(DA_q(x) - DA_q(x-q^{1/3}))$
on $q^{1/3} \N_+$.
Using $A_q(y) - A_q(x) = q^{1/3} \sum_{x \leq z < y} DA_q(z)$, and the fact that
$x \to DA_q(x)$ is non-increasing,  one finds that $|DA_q(x)| \leq x^{-1}$ for all $x$ and $q$.
Together with the identity $\Delta A_q(x) = x A_q(x)$, this shows that both $D A_q(x)$ and
$\Delta A_q$ are bounded over intervals $[0,K]$ uniformly in $q$.
This gives the required equicontinuity to deduce that there is a
continuously differentiable $A:[0,\infty) \to [0,1]$ so that, at least subsequentially,
\[
\max_{x \leq K} |A_q(x) - A(x)| \to 0 \quad \mbox{and} \quad
\max_{x \leq K} |DA_q(x) - A'(x)| \to 0.
\]
We aim to check the limit points $A$ are unique and solve the Airy equation.
Passing to the limit in
\[
\left| q^{1/3} \sum_{x<z \leq y}  z A_q(z) \right| =
\left| q^{1/3} \sum_{x<z \leq y}  \Delta A_q(z) \right|
\leq - DA_q(x) \leq x^{-1}
\]
leads to $\int_x^{\infty}  z A(z) \, dz \leq x^{-1}$ and hence that $A(z) \to 0$ as $z \to \infty$.
The other boundary condition $A(0)=1$ is immediate.
Another passage to the limit in
\[
A_q(x) = 1+ q^{1/3} \sum_{0 \leq y < x} ( D_q(y) + q^{1/2} \sum_{0 \leq z < y} z A_q(z))
\]
gives $A(x) = 1+ \int^x_0 (A'(0) + \int^y_0 x A(z)  dz)  dy$, completing the proof that
$A \in C^2([0,\infty))$, a solution to the Airy equation.  Uniqueness of solutions is
straightforward via a maximum principle argument.

\end{document}